\documentclass[reqno,11pt]{amsart}

\usepackage{amsmath}
\usepackage{amssymb}
\newtheorem{Theorem}{Theorem}[section]
\newtheorem{Lemma}[Theorem]{Lemma}
\newtheorem{Proposition}[Theorem]{Proposition}

\theoremstyle{definition}

\newtheorem{Example}[Theorem]{Example}

\numberwithin{equation}{section}

\begin{document}

\title{Heat kernel expansions on the integers}

\author[F.~A.~Gr\"unbaum]{F.~Alberto~Gr\"unbaum}
\address{F.~A.~Gr\"unbaum, Department of Mathematics, University of 
California, Berkeley, CA 94720-3840, USA}
\email{grunbaum@math.berkeley.edu}

\author[P.~Iliev]{Plamen~Iliev}
\address{P.~Iliev, Department of Mathematics, University of California, 
Berkeley, CA 94720-3840, USA}
\email{iliev@math.berkeley.edu}

\date{December 27, 2001, Revised version: April 21, 2002}

\begin{abstract}
In the case of the heat equation $u_t=u_{xx}+Vu$ on the real line
there are some remarkable potentials $V$ for which the asymptotic
expansion of the fundamental solution becomes a finite sum and gives
an exact formula.

We show that a similar phenomenon holds when one replaces the real line
by the integers. In this case the second derivative is replaced by the 
second difference operator $L_0$. We show if $L$ denotes the result of 
applying a finite number of Darboux transformations to $L_0$ then the 
fundamental solution of $u_t=Lu$ is given by a finite sum of terms 
involving the Bessel function $I$ of imaginary argument.
\end{abstract}

\maketitle

\newcommand{\thref}[1]{Theorem \ref{#1}}
\newcommand{\leref}[1]{Lemma \ref{#1}}
\newcommand{\prref}[1]{Proposition \ref{#1}}
\newcommand{\reref}[1]{Remark \ref{#1}}
\newcommand{\deref}[1]{Definition \ref{#1}}
\newcommand{\exref}[1]{Example \ref{#1}}
\newcommand{\coref}[1]{Corollary \ref{#1}}
\newcommand{\cC}{{\mathcal C}}
\newcommand{\cA}{{\mathcal A}}
\newcommand{\cP}{{\mathcal P}}
\newcommand{\cR}{{\mathcal R}}
\newcommand{\cS}{{\mathcal S}}
\newcommand{\cV}{{\mathcal V}}
\newcommand{\cW}{{\mathcal W}}
\newcommand{\C}{\mathbb C}
\newcommand{\Z}{\mathbb Z}
\newcommand{\Van}{\mathrm{Van}}
\newcommand{\Spec}{\mathrm{Spec}}
\newcommand{\Wr}{\mathrm{Wr}}
\newcommand{\spn}{\mathrm{span}}
\newcommand{\ord}{\mathrm{ord}}
\newcommand{\res}{\mathrm{res}}
\newcommand{\Ker}{\mathrm{Ker}}
\newcommand{\Id}{\mathrm{Id}}
\newcommand{\pd}{\partial}
\newcommand{\La}{\Lambda}
\newcommand{\De}{\Delta}
\newcommand{\expp}{\exp\left(\sum_{i=1}^{\infty}r_iz^i\right)}
\newcommand{\expm}{\exp\left(-\sum_{i=1}^{\infty}r_iz^i \right)}
\newcommand{\Expp}{{\mathrm{Exp}}(n;r,z)}
\newcommand{\Expm}{{\mathrm{Exp}}^{-1}(n;r,z)}


\section{Heat kernel expansions}

The subject of heat kernel expansions on Riemannian manifolds with or without
boundaries cuts across a number of branches of mathematics and serves as an
interesting playground for a whole array of interactions with physics.

It suffices to mention, for instance, the work of M.~Kac \cite{K} on the
issue of recovering the shape of a drum from its pure tones, as well as the
suggestion by H.~P.~McKean and I.~Singer \cite{McKS} that one should be able
to find a ``heat equation proof'' of the index theorem.  For an updated
account, see \cite{BGV,R2}.

In the simpler case of the whole real line, the fundamental solution 
of the equation
\begin{equation*}
u_t = u_{xx} + V(x)u,\quad u(x,0) = \delta_y(x)
\end{equation*}
admits an asymptotic expansion valid for {\em small $t$ and $x$ close to $y$},
in the form
\begin{equation*}
u(x,y,t) \sim \frac {e^{-\frac{(x-y)^2}{4t}}}{\sqrt{4\pi t}} \left( 1 +
\sum_{n=1}^{\infty} H_n(x,y)t^n\right).
\end{equation*}
When $V$ is taken to be a potential such that $L=(d/dx)^2+V$ belongs to
a {\em rank one bispectral ring}, then something remarkable happens, namely 
this 
expansion gives rise to an exact formula consisting of a finite number of
terms and valid {\em{for all $x$, $y$ as well as all $t$}}.  For a few 
examples of this see \cite{G2}, as well as references in \cite{AS}. 
See also \cite{Sch} where the author points to \cite{McKvM} for useful 
connections between the heat equation and KdV invariants. 
For a general discussion of the {\em bispectral problem} see \cite{DG} 
and \cite{W} and for a sample of problems touching upon this area see 
\cite{Be,BK,BV,BW,ChFV,FMTV,GLZ}.

In the case above the potentials $V$'s are the rational solutions of the KdV 
equation decaying at infinity, see \cite{AMcM}. They can be obtained by a
finite number of applications of the Darboux process starting from the 
trivial potential $V=0$, see \cite{AM}. The remarks above rest on the 
basic fact that the Darboux process maps operators of the form 
$d^2/dx^2+ V$ into themselves.

Following some preliminary explorations in \cite{G1} we would like
to use this ``soliton technology'' as a tool for the discovery 
of the general form of a heat kernel expansion on the integers.
In this case, there seems to be no general theory predicting even the 
existence of asymptotic expansions.

If one replaces the real line by the integers and  looks for the fundamental 
solution of
\begin{equation*}
u_t(n,t) = u(n+1,t) - 2u(n,t) + u(n-1,t)\equiv L_0 u(n,t)
\end{equation*}
with $u(n,0) = \delta_{nm}$, one obtains, in terms of the Bessel function 
$I_n(t)$ of imaginary argument, the well-known result
\begin{equation*}
u(n,t) = e^{-2t}I_{n-m}(2t).
\end{equation*}
A nice reference for this is \cite{F} , volume 2. In the context of 
$l^2(\Z)$ it is simplest to study $L_0$ by Fourier methods and obtain 
for the fundamental solution of the heat equation the expression
\begin{equation}                                                 \label{1.1}
\frac{e^{-2t}}{2\pi i}\oint e^{t(x+x^{-1})} x^{n-m}\frac{dx}{x}.
\end{equation}

One can replace the second difference operator $L_0$ above by an
appropriate perturbation of it and look at the corresponding heat equation 
and its fundamental solution. The purpose of this paper is to 
describe the result when the ordinary second difference operator
is subject to a finite number of Darboux factorization steps. In this
case the spectrum is a finite interval and the factorization can be performed 
at each end (see \eqref{2.1}). Formula \eqref{1.1} above will be modified 
properly, in \eqref{5.3}, when $L_0$ is subject to a finite number of 
applications of the Darboux process.

We close this introduction with the simplest nontrivial example and then we 
describe the organization of the paper.

Let us write the operator $L_0$ in the form $L_0=\Lambda-2\Id+\Lambda^{-1}$, 
where $\Lambda$ stands for the customary shift operator, acting on functions 
of a discrete variable $n\in\Z$ by $\Lambda f(n)=f(n+1)$. We can factorize 
$L_0$ as 
\begin{equation}						\label{1.2}
L_0=P_0Q_0,
\end{equation}
where the operators $P_0$ and $Q_0$ are given by
\begin{equation}						\label{1.3}
P_0=\Id-\frac{\tau_{n-1}(\delta)}{\tau_{n}(\delta)}\Lambda^{-1}\text{ and }
Q_0=\Lambda-\frac{\tau_{n+1}(\delta)}{\tau_{n}(\delta)}\Id,
\end{equation}
with $\tau_n(\delta)=n+\delta$. Applying one Darboux step with parameter 
$\delta$ to the operator $L_{0}$ amounts to producing a new operator 
$L_{1,0}$ by exchanging the order of the factors in \eqref{1.2}, i.e.
\begin{equation}						\label{1.4}
L_{1,0}=Q_0P_0=\Lambda
-\Big(2+\frac{1}{\tau_n(\delta)\tau_{n+1}(\delta)}\Big)\Id
+\frac{\tau_{n+1}(\delta)\tau_{n-1}(\delta)}{\tau_n(\delta)^2}\Lambda^{-1}.
\end{equation}
The fundamental solution to
\begin{equation}						\label{1.5}
u_t = L_{1,0}u
\end{equation}
is given by
\begin{equation}						\label{1.6}
u(n,m,t)=\frac{e^{-2t}}{\tau_{m+1}\tau_{n}}
\Big[\tau_m\tau_{n+1}I_{n-m}(2t)-tI_{n-m}(2t)-tI_{n-m+1}(2t)\Big].
\end{equation}
Indeed, it is obvious that $u(n,m,0)=\delta_{nm}$, and equation \eqref{1.5} 
can be verified using well known identities satisfied by the 
Bessel functions (see \eqref{4.2a} and \eqref{4.2b}).

It is important to point out that, when the real line is replaced by 
the integers, the operators obtained from $L_0$ by the Darboux process 
are no longer of the form $L_0$ plus a potential.

The operators $L$, obtained from $L_0$ by a finite number of 
applications of the Darboux process at the ends of the spectrum, have 
been recently determined, see \cite{HI1,HI2}. These operators belong to 
a rank-one commutative ring $\cA_{\cV}$ of difference operators with 
unicursal spectral curve. The common eigenfunction 
$p_n(x)$ to all operators of $\cA_{\cV}$, with spectral parameter $x$, 
is also an eigenfunction to a rank-one commutative ring of differential 
operators in $x$ with solitonic spectral curve, i.e. we have a 
difference-differential bispectral situation. Moreover, the functions 
$p_n(x)$ satisfy an orthogonality relation on the circle. These results will be
summarized in the next section and put to use in later ones.

The ring of operators $\cA_{\cV}$ is also discussed in section 2 and certain 
operators in it are exhibited in section 3. In section 4 we collect some 
properties of Bessel functions and then all of this is used in section 5 
to obtain our main result, \thref{th5.1}. We close the paper with one example 
to illustrate all the steps of the proof of our main result what covers all 
cases when $L$ belongs to a {\em rank one bispectral ring}.

\section{Rank-one bispectral second-order difference operators}

Denote by $\Lambda$ and $\Delta$, respectively, the customary
shift and difference operators, acting on functions
of a discrete variable $n \in \mathbb{Z}$ by
\begin{equation*}
\Lambda f(n)=f(n+1){\text{ and }}\Delta f(n)=f(n+1)-f(n)=(\Lambda-\Id)f(n).
\end{equation*}
The formal adjoint to an operator 
\begin{equation*}
X=\sum_{j}a_j(n)\Lambda^j
\end{equation*}
is defined to be 
\begin{equation*}
X^*=\sum_{j}\Lambda^{-j}\cdot a_j(n)=\sum_{j}a_j(n-j)\Lambda^{-j}.
\end{equation*}

In this section we describe the operator $L_{R,S}$ obtained by successive 
Darboux transformations from the operator $L_0=\Lambda-2\Id+\Lambda^{-1}$.
The symmetric operator $L_0$ admits (as $d^2/dx^2$ did in the case of the 
real line) a unique selfadjoint extension in $l^2(\Z)$ (``limit point'' case 
of Weyl's classification at both end points). Its spectrum is the interval 
$(-4,0)$. The steps of the Darboux transformations are as follows
\begin{align}
&L_{0}=P_{0}Q_{0} \curvearrowright L_{1,0}=Q_{0}P_{0}=P_{1}Q_{1}
\curvearrowright\cdots \curvearrowright L_{R,0}=Q_{R-1}P_{R-1}    \nonumber\\
&\quad L_{R,0}+4\Id=P_{R}Q_{R} \curvearrowright
L_{R,1}+4\Id=Q_{R}P_{R}= P_{R+1}Q_{R+1}
\curvearrowright\cdots
                                                                  \nonumber\\
&\qquad\curvearrowright L_{R,S}+4\Id=Q_{R+S-1}P_{R+S-1}.          \label{2.1}
\end{align}

At each step we perform a lower-upper factorization, as we did in \eqref{1.3}, 
of the corresponding operator and then we produce a new operator by 
interchanging the factors. The operator $L_{R,S}$ is obtained by performing 
$R$ Darboux steps at one end of the spectrum of $L_0$ and 
$S$ steps at the other end. 

For details of the  exposition that  follows we refer the reader to 
\cite{HI2}. 

Let $e(k,\lambda)$ be the linear functional acting on a function 
$g(z)$ by the formula 
\begin{equation*}
\langle e(k,\lambda),g\rangle =g^{(k)}(\lambda),\quad\lambda\in\C,\, k\geq 0.
\end{equation*}
We shall denote by $\Expp$ the exponential function
\begin{equation*}
\Expp=(1+z)^n\expp,
\end{equation*}
where $r=(r_1,r_2,\dots)$. Let us introduce also the functions
\begin{equation}                                                  \label{2.2}
\cS_{j}^{\varepsilon}(n;r)=
        \frac{1}{j!}\langle e(j,\varepsilon-1),\Expp\rangle,
\end{equation}
and write 
\begin{equation}                                                 \label{2.3}
\phi_{j}(n;r)=\cS^1_{2j-1}(n+j-1;r) {\text{ and }} 
\psi_{j}(n;r)=\cS^{-1}_{2j-1}(n+j-1;r).
\end{equation}
When $\varepsilon=1$, the functions $\cS^{1}_{j}(n;r)$ are a shifted version 
of the classical elementary Schur polynomials, defined by 
$\exp\big(\sum_{j=1}^{\infty}r_jz^j\big)=\sum_{j=0}^{\infty}\cS_j(r)z^j$:
\begin{equation*}
\cS^{1}_j(n;r)=\cS_j(r_1+n,r_2-n/2,r_3+n/3,\dots.).
\end{equation*}
For $\varepsilon=-1$, the functions $\cS^{-1}_{j}(n;r)$ are of the form 
\begin{equation*}
\cS^{-1}_j(n;r)=(\text{a polynomial in }n,r_1,r_2,\dots)\,\times
(-1)^n\exp\Bigg(\sum_{j=1}^{\infty}(-2)^jr_j\Bigg).
\end{equation*}
Finally, let us define 
\begin{align}
\tau(n;r)=&\Wr_\Delta(\phi_{1}(n;r),\dots,\phi_{R}(n;r),\psi_{1}(n;r),\dots,
\psi_{S}(n;r))                          \nonumber\\
&\times (-1)^{nS}\exp\left(-S\sum_{i=1}^{\infty}r_i(-2)^i\right),\label{2.4}
\end{align}
where $\Wr_{\Delta}$ denotes the discrete Wronskian with respect to 
the variable $n$:
\begin{equation*}
\Wr_{\Delta}(f_1(n),f_2(n),\dots,f_k(n))
=\det(\Delta^{i-1}f_j(n))_{1\leq i,j\leq k}.
\end{equation*}
The purpose of the exponential factor in \eqref{2.4} is to cancel the 
exponential factor that comes out from the functions $\psi_j(n;r)$, 
$1\leq j\leq S$. With this normalization, the function $\tau(n;r)$ becomes a 
quasipolynomial in the variables $n,r_1,r_2,\dots$, i.e. in general, 
$\tau(n,r)$ depends on infinitely many variables, but there exists a 
positive integer $N$, such that $\tau(n,r)$ is a polynomial in every 
variable of degree at most $N$. 
The operator $L_{R,S}$ in \eqref{2.1} can be expressed in terms of the 
function $\tau(n,r)$ via the formula 
\begin{equation}                                                 \label{2.5}
L_{R,S}=\La+\left(-2+\frac{\partial}{\partial r_1}
\log\frac{\tau(n+1;r)}{\tau(n;r)}\right)\Id+
\frac{\tau(n-1;r)\tau(n+1;r)}{\tau(n;r)^2}\La^{-1}.
\end{equation}
From \eqref{2.4} it follows immediately that $\tau(n,r)$ is an adelic 
tau function\footnote{By an adelic tau function we mean a tau function 
built from a plane belonging to Wilson's adelic Grassmannian \cite{W} via 
the construction in \cite{HI1}.}
 of the $\Delta$KP hierarchy defined by $R+S$ one-point 
conditions, with $R$ conditions at the point zero and $S$ conditions at the 
point $-2$. The reason why the end points of the spectrum of $L_0$, $0$ and 
$-4$, are replaced by $0$ and $-2$ can be traced back to \eqref{2.2} and 
\eqref{2.3}. After a suitable linear change of time variables $r_1,r_2,\dots$ 
the tridiagonal operator $L_{R,S}$ solves the Toda lattice hierarchy, i.e.
\begin{equation*}
\frac{\partial L}{\partial r'_j}=[(L^j)_+,L], \qquad j=1,2,\dots 
\end{equation*}
where $(L^j)_+$ denotes the positive difference part\footnote{Equivalently, 
if we think of $L^j$ as an infinite matrix, $(L^j)_+$ is the upper part of 
this matrix, including the main diagonal.}
of the operator $L^j$, and $\{r'_j\}_{j=1}^{\infty}$ are related to 
$\{r_j\}_{j=1}^{\infty}$ via linear transformation of the form
\begin{equation*}
r_j=r'_j+\sum_{k=j+1}^{\infty}c_{jk}r'_k,
\end{equation*}
see \cite{HI1} for explicit formulas. 
The wave function $w(n;r,z)$ and the adjoint wave function 
$w^*(n;r,z)$ are defined by
\begin{align}
w(n;t,z)&=\frac{\tau(n;r-[z^{-1}])}{\tau(n;r)}\Expp               \nonumber\\
	&=\Bigg(1+\sum_{j=1}^{\infty}w_j(n;r)z^{-j}\Bigg)\Expp,  \label{2.6}\\
\intertext{and}
w^{*}(n;r,z)&=\frac{\tau(n;r+[z^{-1}])}{\tau(n;r)}\Expm           \nonumber\\
	&=\Bigg(1+\sum_{j=1}^{\infty}w^*_j(n;r)z^{-j}\Bigg)\Expm,\label{2.7}
\end{align}
where $[z]=(z,z^2/2,z^3/3,\dots)$. The wave operator $W(n;r)$ and the adjoint 
wave operator $(W^{-1})^*$ can be written in the form
\begin{align}
W(n;r)&=1+\sum_{j=1}^{\infty}w_{j}(n;r)\Delta^{-j}=
Q(\La-\Id)^{-R}(\La+\Id)^{-S},					\label{2.8}\\
(W^{-1})^*(n;r)&=1+\sum_{j=1}^{\infty}w^*_{j}(n+1;r)\Delta^{*\,-j}=
P^*(\La^*-\Id)^{-R}(\La^*+\Id)^{-S},				\label{2.9}
\end{align}
where $P$ and $Q$ are finite-band forward difference operators of order 
$R+S$ satisfying 
\begin{equation}						\label{2.10}
PQ=(\La-\Id)^{2R}(\La+\Id)^{2S},
\end{equation}
see \cite[Theorem 3.3]{HI2}. 

The tau function $\tau(n;r)$ can be viewed as an infinite sequence (indexed
by $n$) of tau functions of the standard KP hierarchy, associated with a flag
of nested subspaces
\begin{equation}                                                  \label{2.11}
\mathcal{V}:\quad \cdots \subset V_{n+1} \subset V_{n} \subset V_{n-1}
    \subset \cdots,
\end{equation}
with
\begin{equation}                                                  \label{2.12}
V_{n}={\mathrm{Span}} \{w(n;0,z),w(n+1;0,z),w(n+2;0,z),\dots\}.
\end{equation}
The plane $V_{0}$ belongs to Wilson's adelic Grassmannian \cite{W}. 
Let $x=z+1$,  
\begin{equation}                                                  \label{2.13}
p_{n}(x)=W(n;r) x^{n}=
w(n,r,x-1)\exp\left(-\sum_{j=1}^{\infty}r_j(x-1)^j\right),
\end{equation}
and similarly 
\begin{equation}                                                  \label{2.14}
p^*_{n}(x)=W^*(n-1;r) x^{-n}=
w^*(n,r,x-1)\exp\left(\sum_{j=1}^{\infty}r_j(x-1)^j\right).
\end{equation}
From \eqref{2.1} it follows easily that $p_n(x)$ are eigenfunctions of the 
operator $L_{R,S}$ with eigenvalue $x-2+x^{-1}$, i.e. 
\begin{equation}                                                  \label{2.15}
L_{R,S}p_n(x)=(x-2+x^{-1})p_n(x).
\end{equation}
Moreover, $p_n(x)$ are also eigenfunctions of a differential operator in 
$x$ (with coefficients independent of $n$) and thus provide a 
difference-differential analog of the rank-one solutions (the KdV family) of 
the bispectral problem considered by Duistermaat and Gr\"unbaum in \cite{DG}. 

Below we discuss the spectral curve and the common eigenfunction of the 
maximal commutative ring of difference operators $\cA_{\cV}$ containing 
$L_{R,S}$. 
The correspondence between commutative rings of difference operators and 
curves was studied by van~Moerbeke and Mumford \cite{vMM} and Krichever 
\cite{Kr} (see also \cite{Mum} where the case of singular curves was 
treated very completely).

Let us introduce the ring $A_{\cV}$ of Laurent polynomials in $x$ that 
preserve the flag $\cV$:
\begin{equation}                                                  \label{2.16}
A_{\cV}=\{f(x)\in\C[x,x^{-1}]:\;f(x)V_n\subset V_{n+k},\;\text{ for some }
k\in\Z\text{ and }\forall n\in\Z\}.
\end{equation}
For each $f(x)\in A_{\cV}$ there exists a finite band operator $L_f$, 
such that 
\begin{equation*}
L_f p_n(x) =f(x) p_n(x).
\end{equation*}
Moreover, if $f(x)=\sum_{j=m_1}^{m_2}a_jx^j$ with $a_{m_1}\neq 0$ and 
$a_{m_2}\neq 0$, then the operator $L_f$ has support $[m_1,m_2]$, i.e.  
\begin{equation*}
L_f=\sum_{j=m_1}^{m_2}b_j(n)\Lambda^j.
\end{equation*}
We can think of $L_f$ as an operator obtained by a Darboux transformation from 
the constant coefficient operator $f(\La)$. The ring $A_{\cV}$ is generated 
by the functions\footnote{Notice that these generators differ slightly from 
the ones chosen in \cite{HI2}.}
\begin{equation}                                                  \label{2.17}
w=\frac{1}{2}\left(x+\frac{1}{x}\right)\text{ and } 
v=\frac{1}{2^{R+S+1}}\frac{(x-1)^{2R+1}(x+1)^{2S+1}}{x^{R+S+1}},
\end{equation}
i.e. $A_{\cV}=\C[w,v]$. Thus the operator $L_{R,S}$ constructed above 
belongs to a maximal rank-one commutative ring of difference operators 
$\cA_{\cV}$ isomorphic to the ring of Laurent polynomials $A_{\cV}$. 
The spectral curve of $\cA_{\cV}$ is
\begin{equation*}
\Spec(\cA_{\cV}):\qquad v^2=(w-1)^{2R+1}(w+1)^{2S+1}.
\end{equation*} 

The operators $Q$ and $P^*$ in formulas \eqref{2.8} and \eqref{2.9}, 
can be expressed in terms of the functions $\{\phi_i(n;r)\}_{i=1}^R$ and 
$\{\psi_j(n;r)\}_{j=1}^S$ as
\begin{align}
Q\,f(n)&=\frac{\Wr_{\Delta}(\phi_{1}(n;r),\dots,\phi_{R}(n;r),\psi_{1}(n;r),
\dots,\psi_{S}(n;r),f(n))}{\Wr_{\Delta}(\phi_{1}(n;r),\dots,\phi_{R}(n;r),
\psi_{1}(n;r),\dots,\psi_{S}(n;r))}                              \label{2.18}\\
\intertext{and}
P^{*}\,f(n)&=\frac{\Wr_{\Delta^{*}}(\phi_{1}^{*}(n;r),\dots,
\phi_{R}^{*}(n;r),\psi_{1}^{*}(n;r),\dots,\psi_{S}^{*}(n;r),f(n))}
{\Wr_{\Delta^{*}}(\phi_{1}^{*}(n;r),\dots,\phi_{R}^{*}(n;r),\psi_{1}^{*}(n;r),
\dots,\psi_{S}^{*}(n;r))},                                       \label{2.19}
\end{align}
with $\phi_{i}^{*}(n;t)=\phi_{i}(n+R+S;t)$, $1 \leq i \leq R$, and 
$\psi_{j}^{*}(n;t)=\psi_{j}(n+R+S;t)$, $1 \leq j \leq S$. Using these
explicit formulas for $Q$ and $P$, one can check that 
\begin{equation}                                                  \label{2.20}
p_n(x^{-1})=\frac{\tau(n+1;r)}{\tau(n;r)}xp^{*}_{n+1}(x).
\end{equation}
Finally, (see \cite[Theorem 5.2]{HI2}) one can prove the following 
orthogonality relation
\begin{equation}                                                  \label{2.21}
\frac{1}{2\pi i} \oint_C p_{n}(x)p_{m}(x^{-1})
\frac{d x}{x}=\frac{\tau(n+1;t)}{\tau(n;t)}\delta_{nm},\; \forall
n,m \in \mathbb{Z},
\end{equation}
where $C$ is any positively oriented simple closed contour surrounding 
the origin, not passing through the points $x=\pm1$. Indeed, $p_n(x)$ are 
rational functions on the Riemann sphere with poles only at 
$x=0,\pm 1,\infty$. The spectral curve $\Spec(\cA_{\cV})$ has cusps at 
$\pm 1$, hence using the fact that the residue of a regular differential 
at a cusp is always zero (see \cite{Serre}), we deduce that 
\begin{equation*}
\res_{x=\pm1}p_n(x)p^*_m(x)dx=0.
\end{equation*}
The proof of \eqref{2.21} now follows from the discrete Kadomtsev-Petviashvili 
bilinear identities and the relation between the wave and the adjoint 
wave functions in \eqref{2.20}.

\section{Certain Laurent polynomials in $A_{\cV}$}

The main result in this section is \prref{pr3.1} below which allows us 
to construct some ``universal'' Laurent polynomials in $A_{\cV}$. This 
result is crucial for the proof of \thref{th5.1}. As we shall see in Section 
5, applying \prref{pr3.1}, we can reduce the infinite sum in the 
computation of the fundamental solution of \eqref{5.1}-\eqref{5.2} to 
a {\it finite\/} sum, modulo some identities among the Bessel functions to 
be discussed in the next section.

\begin{Proposition}\label{pr3.1}                           
Let $T$ be $\geq\max(R,S)$, and let $s_0,s_1,\dots,s_T$ be 
distinct nonzero integers, such that $s_j\equiv s_k\pmod 2$ and 
$s_j+s_k\neq 0$, for $0\leq j,k\leq T$. Then
\begin{equation}                                                  \label{3.1}
F_{s_0,\dots,s_T}(x)=
\sum_{k=0}^T\frac{x^{s_k}}{s_k\prod_{j\neq k}(s_k^2-s_j^2)}
\in A_{\cV}.
\end{equation}
\end{Proposition}

The proof of \prref{pr3.1} is based on two simple lemmas related to the 
Lagrange interpolation polynomial and the Chebyshev polynomials of the 
second kind.

\begin{Lemma}\label{le3.2}                                     
Let $q(n)$ be an odd polynomial in $n$ of degree at most $2T-1$, and 
let $s_0,s_1,\dots,s_T$ be distinct positive numbers. Then 
\begin{equation}                                                  \label{3.2}
\sum_{k=0}^T\frac{q(s_k)}{s_k\prod_{j\neq k}(s_k^2-s_j^2)}=0.
\end{equation}
\end{Lemma}
\begin{proof}[Proof of \leref{le3.2}] 
Consider the Lagrange interpolation polynomial for $q(n)$ at the 
nodes $\pm s_0,\pm s_1,\dots, \pm s_{T}$. Since $\deg q(n)\leq 2T-1$, 
the Lagrange polynomial must be identically equal to $q(n)$, i.e. 
we have 
\begin{equation*}
\sum_{k=0}^Tq(s_k)\frac{n+s_k}{2s_k}
\prod_{j\neq k}\frac{n^2-s_j^2}{s_k^2-s_j^2} +
\sum_{k=0}^Tq(-s_k)\frac{n-s_k}{-2s_k}
\prod_{j\neq k}\frac{n^2-s_j^2}{s_k^2-s_j^2} = q(n).
\end{equation*}
Hence, the coefficient of $n^{2T+1}$ on left-hand side must be zero, which 
gives \eqref{3.2}.
\end{proof}

The Chebyshev polynomials of the second kind are defined by the following 
three term recurrence relation
\begin{equation}                                                  \label{3.3}
2wU_n(w)=U_{n+1}(w)+U_{n-1}(w), \quad n=1,2,\dots,
\end{equation}
with $U_0(w)=1$ and $U_1(w)=2w$. From the last formula one sees immediately 
that $U_n(w)$ is an even/odd polynomial if $n$ is an even/odd positive 
integer, respectively. 

The Chebyshev polynomials of the second kind can be obtained from the 
Jacobi polynomials by taking $\alpha=\beta=\frac{1}{2}$:
\begin{align}
U_n(w)&=(n+1)\,{}_2F_1\Big(
\begin{matrix} -n,n+2\\
\frac{3}{2}\end{matrix}\,\Big\vert\,\frac{1-w}{2}\Big)            \nonumber\\
&=(n+1)\sum_{k=0}^n\frac{2^k}{(2k+1)!}
\Bigg(\prod_{j=1}^k(n+1)^2-j^2)\Bigg)(w-1)^k.                     \label{3.4}
\end{align}
These explicit expressions for $U_n(w)$ will be needed below.
\begin{Lemma}\label{le3.3}                                     
If $s_0,s_1,\dots,s_T$ are distinct positive integers, then the polynomial 
\begin{equation}                                                  \label{3.5}
Q_{s_0,\dots,s_T}(w)=
\sum_{k=0}^{T}\frac{U_{s_k-1}(w)}{s_k\prod_{j\neq k}(s_k^2-s_j^2)}
\end{equation}
is divisible by $(w-1)^T$. Moreover if $s_0\equiv s_1\equiv\cdots\equiv s_T 
\pmod 2$, then 
\begin{equation}                                                  \label{3.6}
(w^2-1)^T/Q_{s_0,\dots,s_T}(w).
\end{equation}
\end{Lemma}
\begin{proof}[Proof of \leref{le3.3}]
From \eqref{3.4} we see that 
\begin{equation*}
U_{n-1}^{(k)}(1)=\frac{2^kk!}{(2k+1)!}n\prod_{j=1}^k(n^2-j^2),
\end{equation*}
i.e. $U_{n-1}^{(k)}(1)$ is an odd polynomial in $n$ of degree $2k+1$. 
From \leref{le3.2} it follows that 
\begin{equation*}
Q_{s_0,\dots,s_T}^{(k)}(1)=0,\quad k=0,1,\dots,T-1,
\end{equation*}
hence $(w-1)^T/Q_{s_0,\dots,s_T}(w)$. If $s_0\equiv s_1\equiv\cdots\equiv s_T 
\pmod 2$, then $Q_{s_0,\dots,s_T}(w)$ is either even or odd polynomial, which 
proves \eqref{3.6}.
\end{proof}

\begin{proof}[Proof of \prref{pr3.1}]
Notice that
\begin{align*}
&F_{s_0,\dots,s_{k-1},s_{k},s_{k+1},\dots,s_T}(x)-
F_{s_0,\dots,s_{k-1},-s_{k},s_{k+1},\dots,s_T}(x)\\
&\qquad=\frac{1}{s_k\prod_{j\neq k}(s_k^2-s_j^2)}(x^{s_k}+x^{-s_k})\in A_{\cV}.
\end{align*}
Thus, without any restriction, we may assume that $s_0,\dots,s_T$ are 
positive integers. We can rewrite the first equation in \eqref{2.17} 
as $x^2=2wx-1$. From this relation, one can easily see by induction that 
\begin{equation}                                                 \label{3.7}
x^k=U_{k-1}(w)x-U_{k-2}(w), \text{ for }k=1,2,\dots,
\end{equation}
with the understanding that $U_{-1}(w)=0$. Using \eqref{3.7}, we can 
write the function in \eqref{3.1} in the form 
\begin{equation}                                                 \label{3.8}
F_{s_0,\dots,s_T}(x)=Q_{s_0,\dots,s_T}(w) x +
\sum_{k=0}^T\frac{U_{s_k-2}(w)}{s_k\prod_{j\neq k}(s_k^2-s_j^2)},
\end{equation}
where $Q_{s_0,\dots,s_T}(w)$ is the polynomial defined by \eqref{3.5}. 
Clearly, the sum in \eqref{3.8} belongs to $A_{\cV}$, so it remains to 
show that
\begin{equation}                                                 \label{3.9}
Q_{s_0,\dots,s_T}(w) x\in A_{\cV}.
\end{equation}
Using \eqref{2.17} we can write $v$ in terms of $w$ and $x$ as
\begin{equation*}
v=(w-1)^R(w+1)^S(x-w),
\end{equation*}
which gives
\begin{equation*}
(w-1)^R(w+1)^Sx=v+w(w-1)^R(w+1)^S\in A_{\cV}.
\end{equation*}
Therefore, if $T\geq\max(R,S)$, we have 
\begin{equation}                                                 \label{3.10}
(w^2-1)^Tx\in A_{\cV}.
\end{equation}
The proof now follows from \leref{le3.3}.
\end{proof}

\section{Some identities satisfied by the Bessel functions}

The Bessel functions of imaginary argument are defined by the generating 
function
\begin{equation}						\label{4.1}
\sum_{k\in\Z}I_k(t)x^k=e^{t(x+x^{-1})/2},
\end{equation}	
whence $I_k(t)=I_{-k}(t)$, and 
$I_k(-t)=(-1)^kI_k(t)$. Differentiating \eqref{4.1} with respect to $x$ and 
$t$, one gets
\begin{subequations}
\begin{align}
kI_k(t)&=\frac{t}{2}\left(I_{k-1}(t)-I_{k+1}(t)\right),         \label{4.2a}\\
\intertext{and}
I'_{k}(t)&=\frac{1}{2}\left(I_{k-1}(t)+I_{k+1}(t)\right),       \label{4.2b}
\end{align}
\end{subequations}
respectively. Similarly, one can show that $I_k(t)$ satisfy the 
modified Bessel equation
\begin{equation}                                                \label{4.3}
\left(t^2\partial_t^2+t\partial_t-(t^2+k^2)\right)I_k(t)=0,
\end{equation}
where $\partial_t=d/dt$.
The main result in this section is \prref{pr4.1} below, which says that if 
$q(j)$ is an odd polynomial of $j$, then the sum
\begin{equation*}
\sum_{\begin{subarray}{c} j> k\\ j\text{ odd/even}\end{subarray}}
q(j)I_j(t),
\end{equation*}
can be written as a {\it finite\/} linear combination of Bessel functions 
of the form 
\begin{equation*}
\sum_{\text{finitely many $j$'s}}\alpha_j(t)I_j(t), 
\end{equation*}
where each coefficient $\alpha_j(t)$ is a polynomial in $t$. 
Let us define a sequence of polynomials $\alpha^n_j(t)$ for $j\in\Z$ and 
$n=0,1,\dots$ as follows
\begin{align}
\alpha^0_j(t)&=\delta_{j0}\,\frac{t}{2}=
\begin{cases} t/2 & \text{if }j=0\\ 
0& \text{if } j\neq 0\end{cases}                                \label{4.4}\\
\intertext{and}
\alpha^{n+1}_j(t)&=(t^2\partial_t^2+t\partial_t)\alpha^n_j(t)+
\Big(t^2\partial_t+\frac{t}{2}\Big)(\alpha^n_{j+1}(t)+\alpha^n_{j-1}(t))+ 
                                                                \nonumber\\
&\qquad\qquad\frac{t^2}{4}(\alpha^n_{j+2}(t)-2\alpha^n_j(t)+\alpha^n_{j-2}(t)).
								\label{4.5}
\end{align}
From the symmetry of the defining relation, it follows that 
$\alpha^n_j(t)=\alpha^n_{-j}(t)$. Moreover, it is clear that 
\begin{equation}                                                \label{4.6}
\alpha^n_j=0 \text{ for }j\notin [-2n,2n], 
\end{equation}
and $\alpha^n_{-2n}(t)=\alpha^n_{2n}(t)=t^{2n+1}/2^{2n+1}$. For arbitrary 
$j\in[-2n,2n]$, $\alpha^n_j(t)$ is a polynomial in $t$ of degree at most 
$2n+1$.

\begin{Proposition}\label{pr4.1}                           
Let $k$ and $n$ be integers, $n\geq 0$. Then
\begin{equation}                                               \label{4.7}
\sum_{\begin{subarray}{c}j>k\\ j\equiv k+1\;({\mathrm{mod}}\; 2)
\end{subarray}}j^{2n+1}I_{j}(t)
=\sum_{s\in\Z}\alpha^n_{s-k}(t)I_s(t)
=\sum_{s=-2n+k}^{2n+k}\alpha^n_{s-k}(t)I_s(t).
\end{equation}
\end{Proposition}
\begin{proof}[Proof by induction] For $n=0$ we have to show that 
\begin{equation}                                               \label{4.8}
\sum_{\begin{subarray}{c}j>k\\ j\equiv k+1\;(\textrm{mod}\; 2) \end{subarray}}
jI_{j}(t)=\frac{t}{2}I_k(t),
\end{equation}
which easily follows from \eqref{4.2a}. Assume that \eqref{4.7} holds for 
some $n$. Then, using the modified Bessel equation \eqref{4.3} we obtain
\begin{align*}
&\sum_{l=0}^{\infty}(k+2l+1)^{2n+3}I_{k+2l+1}(t)\\
&\qquad =(t^2\partial_t^2+t\partial_t-t^2)\,
\sum_{l=0}^{\infty}(k+2l+1)^{2n+1}I_{k+2l+1}(t)\\
&\qquad =(t^2\partial_t^2+t\partial_t-t^2)\,
\sum_{s\in\Z}\alpha^n_{s-k}(t)I_s(t),
\end{align*}
which shows that \eqref{4.7} holds for $n+1$, upon using \eqref{4.2b}.
\end{proof}

\section{Heat kernel expansions on the integers} 
In the case of the real line, the solution of the heat equation is not 
unique unless the class of solutions satisfies a condition of the form
\begin{equation*}
|u(x,t)|\leq c_1 e^{c_2x^2},\quad c_1,c_2\geq 0.
\end{equation*}
When one says that the Gaussian kernel is the fundamental solution of the 
heat equation one (implicitly) assumes that one is considering functions 
of moderate growth at infinity. 
Denote by $u(n,m,t)$ the solution of 
\begin{align}
u_t & = L_{R,S}u,                                                \label{5.1}\\
u|_{t=0}&=\delta_{nm},                                           \label{5.2}
\end{align}
where $L_{R,S}$ is the second-order difference operator (acting on 
functions of a discrete variable $n\in\Z$) constructed in Section 2. 
We make the same implicit assumption here when we say that $u(n,m,t)$ is the 
fundamental solution.

From \eqref{2.15} and \eqref{2.21} and elementary spectral theory, 
it follows that
\begin{equation}                                                 \label{5.3}
u(n,m,t)= \frac{\tau(m)}{\tau(m+1)}
\frac{e^{-2t}}{2\pi i}\oint_C e^{t(x+x^{-1})}p_n(x)p_m(x^{-1})\frac{dx}{x},
\end{equation}
where, for simplicity, we have omitted the dependence on the parameters 
$r=(r_1,r_2,\dots)$. Using \eqref{2.20} we can write also the fundamental 
solution in terms of the wave and the adjoint wave functions as
\begin{equation}                                                 \label{5.4}
u(n,m,t)=\frac{e^{-2t}}{2\pi i}
\oint_C e^{t(x+x^{-1})}p_n(x)p^*_{m+1}(x)\, dx.
\end{equation}

\begin{Theorem}\label{th5.1}                                   
The solution of \eqref{5.1} with initial condition \eqref{5.2} can be 
written in the form
\begin{equation}                                                 \label{5.5}
u(n,m,t)=e^{-2t} \sum_{_{\text{finitely many $j$'s}}}\beta_j(n,m,t)I_j(2t),
\end{equation}
where $\beta_j(n,m,t)$ are polynomials in $t$ of degree at most $2T-1$, with 
$T=\max(R,S)$.
\end{Theorem}

\begin{proof}
From \eqref{2.20} it follows that
\begin{equation*}
u(n,m,t)=\frac{\tau(m)\tau(n+1)}{\tau(m+1)\tau(n)} u(m,n,t).
\end{equation*}
Thus, we can assume that $k=n-m\geq 0$. Around $x=0$, we have the expansion
\begin{equation*}
\frac{p_n(x)p_m(x^{-1})}{x}=\frac{\tau(m+1)}{\tau(m)}p_{n}(x)p_{m+1}^{*}(x)
=\frac{x^{n-m-1}}{\tau(n)\tau(m)}
\big(\tau(n+1)\tau(m)+{\mathrm{O}}(x)\big),
\end{equation*}
which shows that 
\begin{equation}                                                 \label{5.6}
\res_{x=0}\Bigg(x^j p_n(x)p_{m}(x^{-1})\frac{dx}{x}\Bigg)=0
\text{ for } j\geq 1-k.
\end{equation}
For $\epsilon =1,2$ and $i=0,1,\dots,T-1$ denote
\begin{equation}                                                 \label{5.7}
q_{k+\epsilon+2i}(j)=\frac{j}{k+\epsilon+2i}
\prod_{\begin{subarray}{c} l=1\\ l\neq i\end{subarray}}^{T}
\frac{j^2-(k+\epsilon+2l)^2}{(k+\epsilon+2i)^2-(k+\epsilon+2l)^2}.
\end{equation}
$q_{k+\epsilon+2i}(j)$ is an odd polynomial in $j$ of degree $2T-1$. 
We have
\begin{subequations}                                             \label{5.8}
\begin{align}
&e^{t(x+x^{-1})} = \sum_{j\geq 1-k} I_{j}(2t) x^j +I_k(2t)x^{-k} \label{5.8a}\\
&\quad + \sum_{j>k+2T}I_{j}(2t)\Bigg(x^{-j}-\sum_{i=0}^{T-1}
q_{k+\epsilon+2i}(j)x^{-(k+\epsilon+2i)}\Bigg)                   \label{5.8b}\\
&\quad +\sum_{\epsilon=1}^2\sum_{i=0}^{T-1}
\left[\sum_{\begin{subarray}{c} j\geq k+\epsilon+2i\\ j\equiv k+\epsilon 
\;({\mathrm{mod}}\; 2)\end{subarray}} q_{k+\epsilon+2i}(j) 
I_{j}(2t)\right] x^{-(k+\epsilon+2i)},                           \label{5.8c}
\end{align}
\end{subequations}
where in the second sum, for every fixed $j>k+2T$, we choose $\epsilon =1$ or 
$2$, so that $j\equiv k+\epsilon {\pmod 2}$. Denote by $f_j(x)$ the 
Laurent polynomial of $x$ in the second sum in \eqref{5.8}, i.e. 
\begin{equation*}
f_j(x)=x^{-j}-\sum_{i=0}^{T-1}
q_{k+\epsilon+2i}(j)x^{-(k+\epsilon+2i)}.
\end{equation*}
Notice that 
\begin{equation*}
f_j(x)=-j\prod_{l=1}^T(j^2-(k+\epsilon+2l)^2)F_{s_0,s_1,\dots,s_T}(x),
\end{equation*}
where $F_{s_0,s_1,\dots,s_T}(x)$ is the function defined by \eqref{3.1} with 
$s_i=-(k+\epsilon+2i)$, $i=0,1,\dots,T-1$ and $s_T=-j$. Thus, by 
\prref{pr3.1}, $f_j(x)\in A_{\cV}$, and therefore there exists a 
difference operator $L_{f_j}$ with support $[-j,-(k+\epsilon)]$, satisfying
\begin{equation*}
f_j(x)p_n(x)=L_{f_j}p_n(x).
\end{equation*}
Since $n-(k+\epsilon)=m-\epsilon<m$ we see that
\begin{equation}                                                 \label{5.9}
\oint_C f_j(x)p_n(x)p_m(x^{-1})\frac{dx}{x}=0.
\end{equation}
From \eqref{5.3}, \eqref{5.6} and \eqref{5.9} it follows that the infinite 
sums in \eqref{5.8a} and \eqref{5.8b} do not ``contribute'' to the fundamental 
solution $u(n,m,t)$. Finally, the infinite sum in \eqref{5.8c} can be 
rewritten as a finite linear combination of Bessel functions with 
polynomial coefficients, according to \prref{pr4.1}, which completes the 
proof.
\end{proof}

We shall illustrate all steps of the proof by considering the case $R=S=1$.

\begin{Example}\label{ex5.2}                                  
Let $R=S=1$. From \eqref{2.2} and \eqref{2.3} one computes  
\begin{align*}
\phi_1(n;r)&=\cS^1_1(n;r)=n+r_1;\\
\psi_1(n;r)&=\cS^{-1}_1(n;r)=\Bigg(-n+\sum_{j=1}^{\infty}(-2)^{j-1}jr_j\Bigg)
(-1)^n\exp\Bigg(\sum_{j=1}^{\infty}(-2)^jr_j\Bigg).
\end{align*}
Denote for simplicity $\alpha=r_1$ and 
$\beta=\sum_{j=2}^{\infty}(-2)^{j-1}jr_j$. 
The tau function is given by formula \eqref{2.4}:
\begin{equation}                                                 \label{5.10}
\tau(n)=\left\lvert\begin{array}{cc}
n+\alpha   & -n+\alpha+\beta\\
1          & 2n+1-2(\alpha+\beta)
\end{array}\right\rvert.
\end{equation}
The second-order difference operator $L_{1,1}$ is given by formula 
\eqref{2.5}, with $\partial/\partial r_1=\partial/\partial \alpha$.

From \eqref{2.6}, \eqref{2.8}, \eqref{2.13} and \eqref{2.18} we get
\begin{equation}                                                 \label{5.11}
p_n(x)=\frac{x^n}{\tau(n)(x^2-1)}
\left\lvert\begin{array}{ccc}
n+\alpha   & -n+\alpha+\beta      & 1  \\
n+1+\alpha & n+1-\alpha-\beta     & x  \\
n+2+\alpha & -n-2+\alpha+\beta    & x^2
\end{array}\right\rvert.
\end{equation}

From the last formula one can easily deduce that near $x=0$ we can 
expand $p_n(x)p_m(x^{-1})$ as 
\begin{equation}                                                 \label{5.12}
x^{n-m}\Bigg(\frac{\tau(n+1)}{\tau(n)}+
\sum_{j=1}^{\infty}\gamma_j x^j \Bigg),
\end{equation}
where
\begin{equation}                                                 \label{5.13}
\begin{split}
\gamma_j&=-\frac{4}{\tau(n)\tau(m)}\Big[(m-\alpha-\beta+1)(n-\alpha-\beta+1)
(n-m+j)\\
&\qquad\qquad +(-1)^j(m+\alpha+1)(n+\alpha+1)(n-m+j)\Big].
\end{split}
\end{equation}
From \eqref{5.13} and \eqref{4.8} one can see that indeed 
$u(n,m,t)$ is a finite linear combination of Bessel functions. Below, we 
shall illustrate how this can be seen following the proof of 
\thref{th5.1} using just the first few coefficients in \eqref{5.12}.

\prref{pr3.1} tells us that $\forall s,l\neq 0$, such that $s\equiv l\pmod 2$ 
\begin{equation*}
\frac{x^s}{s}-\frac{x^l}{l}\in A_{\cV}.
\end{equation*}
Following \eqref{5.8}, we can write $e^{t(x+x^{-1})}$ as
\begin{subequations}                                           \label{5.14}
\begin{align}
&e^{t(x+x^{-1})} = \sum_{j\geq 1-k} I_{j}(2t)x^j + I_{k}(2t)x^{-k}
                                                               \label{5.14a}\\
&\quad +\sum_{\begin{subarray}{c}j>k+2\\ j\equiv k+1\; ({\mathrm{mod}}\;2) 
\end{subarray}}I_{j}(2t)\Big(x^{-j}-\frac{j}{k+1}x^{-(k+1)}\Big)
                                                               \label{5.14b}\\
&\quad +\sum_{\begin{subarray}{c}j>k+2\\ j\equiv k\; ({\mathrm{mod}}\;2) 
\end{subarray}}I_{j}(2t)\Big(x^{-j}-\frac{j}{k+2}x^{-(k+2)}\Big)
                                                               \label{5.14c}\\
&\quad +\frac{1}{k+1}\left[\sum_{\begin{subarray}{c}j>k\\ j\equiv k+1\; 
({\mathrm{mod}}\;2) \end{subarray}}jI_{j}(2t)\right] x^{-(k+1)}
                                                               \label{5.14d}\\
&\quad +\frac{1}{k+2}\left[\sum_{\begin{subarray}{c}j>k+1\\ j\equiv k\; 
({\mathrm{mod}}\;2) \end{subarray}}jI_{j}(2t)\right] x^{-(k+2)},\label{5.14e}
\end{align}
\end{subequations}
where $k=n-m$. Using \eqref{4.8}, we can write the sums in \eqref{5.14d} 
and \eqref{5.14e} as 
\begin{equation*}
\sum_{\begin{subarray}{c}j>k\\ j\equiv k+1\; 
({\mathrm{mod}}\;2) \end{subarray}}jI_{j}(2t) = t I_k(2t)
\text{ and }
\sum_{\begin{subarray}{c}j>k+1\\ j\equiv k\; 
({\mathrm{mod}}\;2) \end{subarray}}jI_{j}(2t) = t I_{k+1}(2t).
\end{equation*}
Thus \eqref{5.14} can be rewritten as
\begin{subequations}                                           \label{5.15}
\begin{align}
&e^{t(x+x^{-1})} = \sum_{j\geq 1-k} I_{j}(2t)x^j              \label{5.15a}\\
&\quad +\sum_{\begin{subarray}{c}j>k+2\\ j\equiv k+1\; ({\mathrm{mod}}\;2) 
\end{subarray}}I_{j}(2t)\Big(x^{-j}-\frac{j}{k+1}x^{-(k+1)}\Big)
                                                               \label{5.15b}\\
&\quad +\sum_{\begin{subarray}{c}j>k+2\\ j\equiv k\; ({\mathrm{mod}}\;2) 
\end{subarray}}I_{j}(2t)\Big(x^{-j}-\frac{j}{k+2}x^{-(k+2)}\Big)
                                                               \label{5.15c}\\
&\quad + x^{-k}\Big[I_{k}(2t)+\frac{t}{k+1}I_{k}(2t)x^{-1}+
\frac{t}{k+2}I_{k+1}(2t)x^{-2}\Big].                           \label{5.15d}
\end{align}
\end{subequations}
The sums in \eqref{5.15a}, \eqref{5.15b} and \eqref{5.15c} 
``do not contribute'' to the integral \eqref{5.3} (see the proof of 
\thref{th5.1}). Thus 
\begin{equation}                                               \label{5.16}
\begin{split}
&u(n,m,t)=e^{-2t}\res_{x=0}\Big[\Big(I_{k}(2t)+
\frac{t}{k+1}I_{k}(2t)x^{-1}+\frac{t}{k+2}I_{k+1}(2t)x^{-2}\Big)\\
&\qquad\qquad\times\frac{\tau(m)}{\tau(m+1)}\frac{p_n(x)p_m(x^{-1})}{x^{k+1}}
\Big]\\
&\qquad =e^{-2t}\frac{\tau(m)}{\tau(m+1)}
\Big(\frac{\tau(n+1)}{\tau(n)} I_{k}(2t)+\frac{\gamma_1}{k+1}tI_{k}(2t) +
\frac{\gamma_2}{k+2}tI_{k+1}(2t)\Big),
\end{split}
\end{equation}
where $\gamma_1$ and $\gamma_2$ are the coefficients in the expansion 
\eqref{5.12}. Using \eqref{5.13} we get the following explicit formula 
for $u(n,m,t)$
\begin{equation}                                              \label{5.17}
\begin{split}
&u(n,m,t) = \frac{e^{-2t}}{{\tau(m+1)\tau(n)}}
\Big[\tau(m)\tau(n+1)I_{n-m}(2t)\\
&\quad + 4t(\beta+2\alpha)(n+m-\beta+2) I_{n-m}(2t)\\
&\quad - 4t(
2mn-\beta n +2n -\beta m +2m+\beta^2+2\alpha\beta-2\beta+2\alpha^2+2) 
I_{n-m+1}(2t)\Big]
\end{split}
\end{equation}
If we put $\delta=\alpha+1/2$ and let $\beta\rightarrow\infty$
we get 
\begin{equation*}
u(n,m,t)=\frac{e^{-2t}}{\tau_{m+1}\tau_{n}}
\Big[\tau_m\tau_{n+1}I_{n-m}(2t)-tI_{n-m}(2t)-tI_{n-m+1}(2t)\Big],
\end{equation*}
where $\tau_n=n+\delta$. Notice that this is exactly formula \eqref{1.6}
for the fundamental solution computed in the introduction ($R=1$, $S=0$).
\end{Example}

The referee has raised the interesting possibility of using the formula
\begin{equation*}
e^{tAB}=1+A\left(\int_0^te^{sBA}\,ds\right)B
\end{equation*}
to get an alternative proof of \thref{th5.1}. This remains a challenging 
problem.

\bigskip
{\bf Acknowledgments.} We thank Henry P.~McKean and the referee for 
suggestions that led to an improved version of this paper.



\begin{thebibliography}{xx}

\bibitem{AM} M.~Adler and J.~Moser, {\em On a class of polynomials connected 
with the Korteweg-de Vries equation}, Comm. Math. Phys. {\bf 61} (1978), 1--30.

\bibitem{AMcM}  H.~Airault, H.~P.~McKean  and J.~Moser, 
{\em Rational and elliptic solutions of the Korteweg-de Vries equation and 
a related many-body problem}, Comm. Pure Appl. Math. {\bf 30} (1977), 95--148.

\bibitem{AS} I.~G.~Avramidi and R.~Schimming, {\em Heat kernel
coefficients for the matrix Schr\"odinger operator}, J. Math. Phys. 
{\bf 36} (1995), 5042--5054.

\bibitem{Be} Y.~ Berest, Huygens principle and the bispectral
problem, In: {\em The bispectral problem (Montreal, PQ, 1997)}, 11--30, 
CRM Proc. Lecture Notes, {\bf 14}, Amer. Math. Soc., Providence, RI, 1998.

\bibitem{BK} Y.~Berest and A.~Kasman, {\em $\mathcal{D}$-modules and Darboux
transformations}, Lett. Math. Phys. {\bf 43} (1998), 279--294.

\bibitem{BV} Y.~Berest and A.~P.~Veselov, {\em The Huygens principle and
integrability}, (Russian) Uspekhi Mat. Nauk {\bf 49} (1994), 7--78, 
{\it translation in} Russian Math. Surveys {\bf 49} (1994), 5--77.

\bibitem{BW} Y.~Berest and G.~Wilson, {\em Classification of rings of
differential operators on affine curves}, Internat. Math. Res. Notices 
{\bf 2} (1999), 105--109.

\bibitem{BGV} N.~Berline, E.~Getzler and M.~Vergne, {\em Heat
kernels and Dirac operators}, Grundlehren der Mathematischen Wissenschaften, 
{\bf 298}, Springer-Verlag, Berlin, 1992.

\bibitem{ChFV} O.~A.~Chalykh, M.~V.~Feigin and A.~P.~Veselov, {\em
Multidimensional Baker--Akhiezer functions and Huygens' principle}, 
Comm. Math. Phys. {\bf 206} (1999), 533--566.

\bibitem{DG} J.~J.~Duistermaat and F.~A.~Gr\"unbaum, {\em Differential
equations in the spectral parameter}, Comm. Math. Phys. {\bf 103} (1986),
177--240.

\bibitem{FMTV} G.~Felder, Y.~Markov, V.~Tarasov and A.~Varchenko,
{\em Differential equations compatible with KZ equations}, 
Math. Phys. Anal. Geom. {\bf 3} (2000), 139--177.

\bibitem{F} W.~Feller, {\em An Introduction to Probability Theory and Its 
Applications}, Volume 2, J. Wiley, 1990.

\bibitem{GLZ} Ya.~I.~Granovskii, I.~M.~Lutzenko and A.~S.~Zhedanov,
{\em Mutual integrability, quadratic algebras, and dynamical symmetry}, 
Ann. Physics {\bf 217} (1992), 1--20.

\bibitem{G1} F.~A.~Gr\"unbaum, The bispectral problem: an overview, In: 
{\em Special Functions 2000: Current Perspective and Future Directions}, 
Eds. J.~Bustoz et al. (2001), 129--140.
\bibitem{G2} F.~A.~Gr\"unbaum, Some bispectral musings, In: 
{\em The bispectral problem (Montreal, PQ, 1997)}, 11--30, 
CRM Proc. Lecture Notes, {\bf 14}, Amer. Math. Soc., Providence, RI, 1998.

\bibitem{HI1} L.~Haine and P.~Iliev, {\em Commutative rings of difference 
operators and an adelic flag manifold}, Internat. Math. Res. Notices {\bf 6} 
(2000), 281--323.
\bibitem{HI2} L.~Haine and P.~Iliev, {\em A rational analogue of the 
Krall polynomials}, Kowalevski Workshop on Mathematical Methods of
Regular Dynamics (Leeds, 2000), J. Phys. A: Math. Gen. {\bf 34} (2001), 
2445--2457.

\bibitem{K} M.~Kac, {\em Can one hear the shape of a drum?}, Amer. Math.
Monthly {\bf 73} (1966), 1--23.

\bibitem{Kr} I.~M.~Krichever, {\em Algebraic curves and non-linear 
difference equations}, (Russian) Uspekhi Mat. Nauk {\bf 33} (1978), 215--216, 
{\it translation in} Russ. Math. Surveys {\bf 33} (1978), 255-256.

\bibitem{McKS} H.~P.~McKean and I.~Singer, {\em Curvature and the
eigenvalues of the Laplacian}, J. Differential Geometry {\bf 1} (1967),
43--69.

\bibitem{McKvM} H.~P.~McKean and P.~van~Moerbeke, {\em The spectrum of Hill's 
equation}, Invent. Math. {\bf 30} (1975), 217--274.

\bibitem{Mum} D.~Mumford, An algebro-geometric construction of 
commuting operators and of solutions to the Toda lattice equation, 
Korteweg-de Vries equation and related non-linear equations, 
in: M.~Nagata (ed.), {\em Proceedings of International Symposium on Algebraic 
Geometry (Kyoto 1977)}, Kinokuniya Book Store, Tokyo, 1978, 115--153.

\bibitem{R2} S.~Rosenberg, {\em The Laplacian on a Riemannian Manifold. 
An introduction to analysis on manifolds}, 
London Math. Society Student Texts {\bf 31}, Cambridge University Press, 
Cambridge, 1997.

\bibitem{Sch} R.~Schimming {\em An explicit expression for the Korteweg-de 
Vries hierarchy}, Z. Anal. Anwendungen {\bf 7} (1988), 203-214.

\bibitem{Serre} J.-P. Serre, {\em Groupes alg\'ebriques et corps de
classes}, Hermann, 1959.

\bibitem{vMM} P.~van~Moerbeke and D.~Mumford, {\em The spectrum of 
difference operators and algebraic curves}, Acta Math. {\bf 143} (1979), 
93--154.

\bibitem{W} G.~Wilson, {\em Bispectral commutative ordinary differential
operators}, J. Reine Angew. Math. {\bf 442} (1993), 177--204.

\end{thebibliography}
\end{document}